\documentclass[12pt]{article}
\usepackage{amsmath,amsfonts,amssymb,amsthm}
\begin{document}

\newcommand{\1}{{{\bf 1}}}
\newcommand{\id}{{\rm id}}
\newcommand{\Hom}{{\rm Hom}\,}
\newcommand{\End}{{\rm End}\,}
\newcommand{\Res}{{\rm Res}\,}
\newcommand{\Image}{{\rm Im}\,}
\newcommand{\Ind}{{\rm Ind}\,}
\newcommand{\Aut}{{\rm Aut}\,}
\newcommand{\Ker}{{\rm Ker}\,}
\newcommand{\gr}{{\rm gr}}
\newcommand{\Der}{{\rm Der}\,}

\newcommand{\Z}{\mathbb{Z}}
\newcommand{\Q}{\mathbb{Q}}
\newcommand{\C}{\mathbb{C}}
\newcommand{\N}{\mathbb{N}}

\newcommand{\g}{\mathfrak{g}}
\newcommand{\h}{\mathfrak{h}}
\newcommand{\wt}{{\rm wt}\;}
\newcommand{\CR}{\mathcal{R}}
\newcommand{\D}{\mathcal{D}}
\newcommand{\E}{\mathcal{E}}
\newcommand{\F}{\mathcal{F}}
\newcommand{\Lie}{\mathcal{L}}
\newcommand{\z}{\bf{z}}
\newcommand{\bflam}{\bf{\lambda}}

\def \<{\langle} 
\def \>{\rangle}
\def \be{\begin{equation}\label}
\def \ee{\end{equation}}
\def \bex{\begin{exa}\label}
\def \eex{\end{exa}}
\def \bl{\begin{lem}\label}
\def \el{\end{lem}}
\def \bt{\begin{thm}\label}
\def \et{\end{thm}}
\def \bp{\begin{prop}\label}
\def \ep{\end{prop}}
\def \br{\begin{rem}\label}
\def \er{\end{rem}}
\def \bc{\begin{coro}\label}
\def \ec{\end{coro}}
\def \bd{\begin{de}\label}
\def \ed{\end{de}}

\newtheorem{thm}{Theorem}[section]
\newtheorem{prop}[thm]{Proposition}
\newtheorem{coro}[thm]{Corollary}
\newtheorem{conj}[thm]{Conjecture}
\newtheorem{exa}[thm]{Example}
\newtheorem{lem}[thm]{Lemma}
\newtheorem{rem}[thm]{Remark}
\newtheorem{de}[thm]{Definition}
\newtheorem{hy}[thm]{Hypothesis}
\makeatletter
\@addtoreset{equation}{section}
\def\theequation{\thesection.\arabic{equation}}
\makeatother
\makeatletter

\begin{flushright}
\today
\\
\end{flushright}

\begin{center}{\Large \bf Abelianizing vertex algebras}\\
\vspace{0.5cm}
Haisheng Li\footnote{Partially supported by an NSA grant}\\
Department of Mathematical Sciences, Rutgers University, Camden, NJ 08102\\
and\\
Department of Mathematics, Harbin Normal University, Harbin, China
\end{center}

\begin{abstract}
To every vertex algebra $V$ we associate a canonical decreasing
sequence of subspaces and prove that the associated graded
vector space $\gr (V)$ is naturally a vertex Poisson algebra, in
particular a commutative vertex algebra.
We establish a relation between this decreasing sequence and the sequence
$C_{n}$ introduced by Zhu.  By using the (classical) algebra $\gr(V)$,
we prove that for any vertex algebra $V$, $C_{2}$-cofiniteness 
implies $C_{n}$-cofiniteness for all $n\ge 2$.
We further use $\gr(V)$ to study generating subspaces of certain types 
for lower truncated $\Z$-graded vertex algebras.
\end{abstract}

\baselineskip=16pt

\section{Introduction}
Just as with classical (associative or Lie) algebras,
abelian or commutative vertex algebras (should be) are 
the simplest objects in the category of vertex algebras.
It was known (see [B]) that commutative vertex algebras
exactly amount to differential algebras, namely
unital commutative associative algebras equipped with a derivation.
Related to the notion of commutative vertex algebra,
is the notion of vertex Poisson algebra (see [FB]),
where a vertex Poisson algebra structure
combines a commutative vertex algebra structure, or equivalently, a
differential algebra structure, with a vertex Lie algebra structure
(see [K], [P]). As it was showed in [FB], vertex Poisson algebras 
can be considered as classical limits of vertex algebras.

In the classical theory, a well known method to abelianize an
associative algebra is to use a good increasing filtration and then
consider the associated graded vector space. A typical example is
about the universal enveloping algebra $U(\g)$
of a Lie algebra $\g$ with the filtration $\{U_{n}\}$,
 where for $n\ge 0$, $U_{n}$ is linearly spanned by
the vectors $a_{1}\cdots a_{m}$ for $m\le n,\; a_{1},\dots,a_{m}\in \g$.
In this case, the associated graded Poisson algebra $\gr U(\g)$
is naturally a Poisson algebra and 
the well known Poincar\'{e}-Birkhoff-Witt
theorem says that the associated graded Poisson algebra $\gr U(\g)$
is canonically isomorphic to the symmetric algebra $S(\g)$ 
which is also a Poisson algebra. This
result and the canonical isomorphism have played 
a very important role in Lie theory.

Motivated by this classical result, 
in \cite{li-vpa} we introduced and studied a notion of what we called
good increasing filtration for a vertex algebra $V$ and
we proved that the associated graded vector space $\gr V$ of $V$ 
with respect to a good increasing filtration
is naturally a vertex Poisson algebra.
Furthermore, for any $\N$-graded vertex algebra $V=\coprod_{n\in
\N}V_{(n)}$ with $V_{(0)}=\C {\bf 1}$, we constructed a canonical good
increasing filtration of $V$.  
This increasing filtration was essentially used 
in [KL], [NG], [Bu1,2], [ABD] and [NT] in the study on
generating subspaces of $V$ with a certain property 
analogous to the well known
Poincar\'{e}-Birkhoff-Witt spanning property.

In this paper, we introduce and study ``good'' {\em decreasing}
filtrations for vertex algebras.  To any vertex algebra $V$ we
associate a canonical decreasing sequence $\E$ of subspaces 
$E_{n}$ for $n\ge 0$
and we prove that the associated graded vector space $\gr_{\E}(V)$
is naturally an
$\N$-graded vertex Poisson algebra, where for $n\in \Z$,
$E_{n}$ is linearly spanned by the vectors
$$u^{(1)}_{-1-k_{1}}\cdots u^{(r)}_{-1-k_{r}}v$$
for $r\ge 1,\; u^{(i)},v\in V,\;k_{i}\ge 0$
with $k_{1}+\cdots +k_{r}\ge n$.
Notice that unlike the increasing filtration
which uses the weight grading,
this decreasing sequence uses only the
vertex algebra structure.

For any vertex algebra $V$, there has been 
a fairly well known decreasing sequence ${\mathcal{C}}=\{C_{n}\}_{n\ge 2}$
introduced by Zhu [Z1,2], where for $n\ge 2$, $C_{n}$ is linearly
spanned by the vectors $u_{-n}v$ for $u,v\in V$.
The notion of $C_{2}$ was introduced and used in the fundamental 
study of Zhu on modular invariance, where the finiteness of
$\dim V/C_{2}$ played a crucial role.
It was showed in [Z2] that $V/C_{2}$ has a natural
Poisson algebra structure.  

In this paper, we relate our decreasing
sequence $\E$ with Zhu's sequence ${\mathcal{C}}$.  
In particular,
we show that $C_{2}=E_{1}$ and $C_{3}=E_{2}$.
We then show that the degree zero subspace 
$E_{0}/E_{1}$ of $\gr_{\E}(V)$, which is naturally 
a Poisson algebra, is exactly the
Zhu's Poisson algebra $V/C_{2}$. We further show that
$\gr_{\E}(V)$ as a differential algebra is generated by
the degree zero subalgebra $V/C_{2}$. 
As an application, we
show that for {\em any} vertex algebra $V$, if $V$ is
$C_{2}$-cofinite, then $V$ is $E_{n}$-cofinite and 
$C_{n+2}$-cofinite for all $n\ge 0$. Similarly we show that
if $V$ is a $C_{2}$-cofinite vertex algebra and if
$W$ is a $C_{2}$-cofinite $V$-module, then $W$ is $C_{n}$-cofinite for
all $n\ge 2$.

Under the assumption that
$V$ is an $\N$-graded vertex algebra with $\dim V_{(0)}=1$, it
has been proved before by [GN] (see also [NT], [Bu1,2])
that $C_{2}$-cofiniteness implies $C_{n+2}$-cofinite for all $n\ge 0$. 
On the other hand, the original method of [GN]
and [KL] used this assumption in an essential way.

As we show in this paper, for certain vertex algebras,
both sequences $\E$ and ${\mathcal{C}}$ are trivial 
in the sense that $E_{n}=C_{n+2}=V$ for all
$n\ge 0$. On the other hand, by using the connection between
the two decreasing sequences
we prove that if $V=\coprod_{n\ge t}V_{(n)}$ 
is a lower truncated $\Z$-graded vertex algebra such as a vertex
operator algebra in the sense of [FLM] and [FHL], 
then for any $k$, $C_{n}, E_{n}\subset \coprod_{m\ge k}V_{(m)}$
for $n$ sufficiently large.
Consequently, $\cap_{n\ge 0}E_{n}=\cap_{n\ge 0}C_{n+2}=0$.
(In this case, both sequences are filtrations.) 
Furthermore, using this result and $\gr_{\E}(V)$ 
we show that if a graded subspace
$U$ of $V$ gives rises to a generating subspace of $V/C_{2}$ as an
algebra, then $U$ generates $V$ with a certain spanning
property. Similar results have been obtained before
in [KL], [NG], [Bu1,2] and [NT] under a stronger condition.

This paper is organized as follows: In Section 2, we define the
sequence $\E$ and show that the associated graded vector space is an
$\N$-graded vertex Poisson algebra. In Section 3, we relate the
sequences $\E$ and ${\mathcal{C}}$. In Section 4, we study 
generating subspaces of certain types 
for lower truncated $\Z$-graded vertex algebras.

\section{Decreasing sequence $\E$ and
the vertex Poisson algebra $\gr_{\E}(V)$}

In this section we first recall the definition of a vertex Poisson algebra
from [FB] and we then construct a canonical
decreasing sequence $\E$ for each vertex algebra $V$ and show that the
associated graded vector space $\gr_{\E}V$ is naturally a vertex
Poisson algebra.  We also show that if $V$ is an $\N$-graded vertex
algebra, then the sequence $\E$ is indeed a filtration of $V$.

Let $V$ be a vertex algebra.
We have Borcherds' commutator formula and iterate formula:
\begin{eqnarray}
[u_{m},v_{n}]&=&\sum_{i\ge 0}\binom{m}{i}(u_{i}v)_{m+n-i},
\label{ecommutator}\\
(u_{m}v)_{n}w&=&\sum_{i\ge 0}(-1)^{i}\binom{m}{i}
\left( u_{m-i}v_{n+i}w-(-1)^{m}v_{m+n-i}u_{i}w\right)\label{eiterate}
\end{eqnarray}
for $u,v,w\in V,\; m,n\in \Z$.
Define a (canonical) linear operator $\mathcal{D}$ on $V$ by
\begin{eqnarray}
{\mathcal{D}}(v)=v_{-2}\1\;\;\;\mbox{ for }v\in V.
\end{eqnarray}
Then
\begin{eqnarray}
Y(v,x){\bf 1}=e^{x\D}v\;\;\;\mbox{ for }v\in V.
\end{eqnarray}
Furthermore,
\begin{eqnarray}\label{eDproperty}
[\D,v_{n}]=(\D v)_{n}=-n v_{n-1}
\end{eqnarray}
for $v\in V,\; n\in \Z$. (See for example [LL] for an exposition 
of such facts.)

A vertex algebra $V$ is called a
{\em commutative vertex algebra}if
\begin{eqnarray}\label{ecommutative-va}
[u_{m},v_{n}]=0\;\;\;\mbox{ for }u,v\in V,\; m,n\in \Z.
\end{eqnarray}
It is well known (see [B], [FHL]) that (\ref{ecommutative-va}) 
is equivalent to that
\begin{eqnarray}
u_{n}=0\;\;\;\mbox{ for }u\in V,\; n\ge 0.
\end{eqnarray}

\br{rborcherds-con} 
{\em Let $A$ be any unital commutative associative algebra
with a derivation $d$. Then one has a commutative vertex algebra 
structure on $A$ with $Y(a,x)b=(e^{xd}a)b$ for $a,b\in A$ 
and with the identity $1$ as the vacuum vector (see [B]). 
On the other hand, let $V$ be any commutative vertex algebra. Then
$V$ is naturally a commutative associative algebra with $u\cdot
v=u_{-1}v$ for $u,v\in V$ and with ${\bf 1}$ as the identity and with
$\D$ as a derivation.  Furthermore, $Y(u,x)v=(e^{x\D}u)v$ for $u,v\in
V$. Therefore, a
commutative vertex algebra exactly amounts to a unital commutative
associative algebra equipped with a derivation, which is often called
a {\em differential algebra}.}  
\er

A vertex algebra $V$ equipped with a $\Z$-grading
$V=\coprod_{n\in \Z}V_{(n)}$ is called a {\em $\Z$-graded vertex algebra}
if ${\bf 1}\in V_{(0)}$ and if for $u\in V_{(k)}$ with $k\in \Z$ 
and for $m,n\in \Z$,
\begin{eqnarray}
u_{m}V_{(n)}\subset V_{(n+k-m-1)}.
\end{eqnarray}
We say that a $\Z$-graded vertex algebra $V=\coprod_{n\in \Z}V_{(n)}$ is 
{\em lower truncated} if $V_{(n)}=0$ for $n$ sufficiently small.
In particular, every vertex operator algebra in the sense of
[FLM] and [FHL] is a lower truncated $\Z$-graded vertex algebra.
An {\em $\N$-graded vertex algebra} is defined in the obvious way.
We say that a vertex algebra $V$ is {\em $\Z$-gradable} ({\em $\N$-gradable}) if
there exists a $\Z$-grading ($\N$-grading) such that $V$ becomes 
$\Z$-graded ($\N$-graded) vertex algebra.
We see that a commutative $\Z$-graded vertex algebra is naturally
a $\Z$-graded differential algebra.

The following definition of the notion of vertex Lie algebra is due
to [K] and [P]:

\bd{dvla}
{\em A {\em vertex Lie algebra} is a vector space $V$ 
equipped with a linear operator $D$ and a linear map
\begin{eqnarray}
Y_{-}: V&\rightarrow& \Hom (V,x^{-1}V[x^{-1}]),\nonumber\\
v&\mapsto& Y_{-}(v,x)=\sum_{n\ge 0}v_{n}x^{-n-1}
\end{eqnarray}
such that for $u,v\in V,\; m,n\in \N$,
\begin{eqnarray}
& &(Dv)_{n}=-nv_{n-1},\label{evla1}\\
& &u_{m}v=\sum_{i=0}^{m}(-1)^{m+i+1}\frac{1}{i!}D^{i}v_{m+i}u,\\
& &[u_{m},v_{n}]=\sum_{i=0}^{m}\binom{m}{i} (u_{i}v)_{m+n-i}.\label{evla3}
\end{eqnarray}}
\ed

A {\em module} (see [K]) for a vertex Lie algebra $V$ is a vector space $W$ 
equipped with a linear map
\begin{eqnarray}
Y_{-}^{W}: V&\rightarrow& \Hom (W,x^{-1}W[x^{-1}]),\nonumber\\
v&\mapsto& Y_{-}^{W}(v,x)=\sum_{n\ge 0}v_{n}x^{-n-1}
\end{eqnarray}
such that (\ref{evla1}) and (\ref{evla3}) hold.

Recall the following notion of vertex Poisson algebra from [FB] (cf. [DLM]):

\bd{dvpa}
{\em A {\em vertex Poisson algebra} 
is a commutative vertex algebra $A$,
or equivalently, a (unital) commutative associative algebra
equipped with a derivation $\partial$, equipped with
a vertex Lie algebra structure $(Y_{-},\partial)$ such that
\begin{eqnarray}\label{evpaLie}
Y_{-}(a,x)\in x^{-1}(\Der A)[[x^{-1}]]\;\;\;\mbox{ for }a\in A.
\end{eqnarray}}
\ed

A {\em module} for a vertex Poisson algebra $A$ is a vector space $W$
equipped with a module structure for $A$ as an associative algebra and 
a module structure for $A$ as a vertex Lie algebra such that
\begin{eqnarray}
Y_{-}^{W}(u,x)(vw)=(Y_{-}^{W}(u,x)v)w+vY_{-}^{W}(u,x)w
\end{eqnarray}
for $u,v\in V,\; w\in W$.

The following result obtained in \cite{li-vpa} gives a construction of
vertex Poisson algebras from vertex algebras through certain
increasing filtrations:

\bp{pold-inc}
Let $V$ be a vertex algebra and let $\E=\{E_{n}\}_{n\in \Z}$ be a good increasing filtration
of $V$ in the sense that ${\bf 1}\in E_{0}$,
\begin{eqnarray}
u_{n}E_{s}\subset E_{r+s}
\end{eqnarray}
for $u\in E_{r},\; r,s,n\in \Z$ and
\begin{eqnarray}
u_{n}E_{s}\subset E_{r+s-1}\;\;\;\mbox{ for }n\ge 0.  
\end{eqnarray}
Then the associated graded vector space $\gr_{\E}V=\coprod_{n\in \Z}E_{n+1}/E_{n}$ 
is naturally a vertex Poisson algebra with
\begin{eqnarray}
& & (u+E_{m-1})(v+E_{n-1})=u_{-1}v+E_{m+n-1},\\
& &\partial (u+E_{m-1})=\D u+E_{m-1},\\
& &Y_{-}(u+E_{m-1})(v+E_{n-1})=\sum_{r\ge 0}(u_{r}v+E_{m+n-2})x^{-r-1}
\end{eqnarray}
for $u\in E_{m},\; v\in E_{n}$ with $m,n\in \Z$.
\ep

Furthermore, the following construction of
good increasing filtrations was also given in \cite{li-vpa}:

\bt{told-inc}
Let $V=\coprod_{n\in \N}V_{(n)}$ be an $\N$-graded vertex algebra
such that $V_{(0)}=\C {\bf 1}$. Let $U$ be a graded subspace of 
$V_{+}=\coprod_{n\ge 1}V_{(n)}$
such that 
$$V={\rm span}\{u^{(1)}_{-k_{1}}\cdots u^{(r)}_{-k_{r}}{\bf 1}
\;|\; r\ge 0,\; u^{(i)}\in U,\; k_{i}\ge 1\}.$$
In particular, we can take $U=V_{+}$.
For any $n\ge 0$, denote by $E_{n}^{U}$ the subspace of $V$ 
linearly spanned by the vectors
$$u^{(1)}_{-k_{1}}\cdots u^{(r)}_{-k_{r}}{\bf 1}$$
for $r\ge 0$, for homogeneous vectors $u^{(1)},\dots,u^{(r)}\in U$ 
and for $k_{1},\dots,k_{r}\ge 1$ with $\wt u^{(1)}+\cdots +\wt u^{(r)}\le n$.
Then the sequence $\E_{U}=\{ E_{n}^{U}\}$ is 
a good increasing filtration of $V$.
Furthermore, $\E_{U}$ does not depend on $U$.
\et

Next, we give a construction of vertex Poisson algebras from
vertex algebras using {\em decreasing} filtrations.
First, we formulate the following general result, which is
similar to Proposition \ref{pold-inc} and  which is classical in nature:

\bp{pabstract}
Let $V$ be any vertex algebra and let $\F=\{ F_{n}\}_{n\ge 0}$ be a decreasing sequence of 
subspaces of $V$ such that ${\bf 1}\in F_{0}$ and 
\begin{eqnarray}\label{edfil-condition}
u_{n}v\in F_{r+s-n-1}\;\;\;\mbox{ for }u\in F_{r},\; v\in F_{s},\;r,s\in \N,\;\; n\in \Z,
\end{eqnarray}
where by convention $F_{m}=V$ for $m<0$.
Then the associated graded vector space $\gr_{\F}V=\coprod_{n\ge 0}F_{n}/F_{n+1}$
is naturally an $\N$-graded vertex algebra with
\begin{eqnarray}\label{enva-def}
(u+F_{r+1})_{n}(v+F_{s+1})=u_{n}v+F_{r+s-n}
\end{eqnarray}
for $u\in F_{r},\; v\in F_{s},\; r,s\in \N,\; n\in \Z$ 
and with ${\bf 1}+F_{1}\in F_{0}/F_{1}$ as the vacuum vector. 
Furthermore, $\gr_{\F}V$ is commutative if and only if 
\begin{eqnarray}\label{evpacondition}
u_{n}v\in F_{r+s-n}\;\;\;\mbox{ for }u\in F_{r},\; v\in F_{s},\;r,s,n\in \N.
\end{eqnarray}
Assume (\ref{edfil-condition}) and (\ref{evpacondition}). 
Then the commutative vertex algebra $\gr_{\F}V$ 
is a vertex Poisson algebra where
\begin{eqnarray}\label{evla-def}
& &\partial (u+F_{r+1})=\D u +F_{r+2},\\
& &Y_{-}(u+F_{r+1},x)(v+F_{s+1})=\sum_{n\ge 0}(u_{n}v+F_{r+s-n+1})x^{-n-1}
\end{eqnarray}
for $u\in F_{r},\; v\in F_{s}$ with $r,s\in \N$.
\ep

\begin{proof} Notice that the condition (\ref{edfil-condition}) guarantees that
the operations given in (\ref{enva-def}) are well defined.
Just as with any classical algebras, 
it is straightforward to check that $\gr_{\F}V$ is an $\N$-graded vertex algebra and
it is also clear that $\gr_{\F}V$ is commutative if and only if (\ref{evpacondition}) holds.
Assuming (\ref{edfil-condition}) and (\ref{evpacondition}) we have 
a commutative associative $\N$-graded algebra $\gr_{\F}(V)$
with derivation $\partial$ defined by
$$\partial (u+F_{n+1})=(u+F_{n+1})_{-2}({\bf 1}+F_{1})=u_{-2}{\bf 1}+F_{n+2}=\D u+F_{n+2},$$
noticing that by (\ref{edfil-condition}) we have
$\D u=u_{-2}{\bf 1}\in F_{n+1}$. 
The condition  (\ref{evpacondition}) guarantees that the linear map $Y_{-}$
in (\ref{evla-def}) is well defined. It is straightforward to check that $\gr_{\F}(V)$
equipped with $Y_{-}$ and $\partial$ is a vertex Lie algebra.

Now, we check the compatibility condition (\ref{evpaLie}).
Let $u\in F_{r},\; v\in F_{s},\; w\in F_{k}$ with $r,s,k\in \N$. 
For $m\ge 0$, using the Borcherds' commutator formula for $V$ we have
\begin{eqnarray}
u_{m}(v_{-1}w)
&=&v_{-1}(u_{m}w)+\sum_{i\ge 0}\binom{m}{i}(u_{i}v)_{m-1-i}w\nonumber\\
&=&v_{-1}(u_{m}w)+(u_{m}v)_{-1}w+\sum_{i=0}^{m-1}\binom{m}{i}(u_{i}v)_{m-1-i}w.
\end{eqnarray}
For $0\le i\le m-1$, using (\ref{evpacondition}) (twice) we have
$$(u_{i}v)_{m-1-i}w\in F_{r+s+k-m+1}.$$
Thus
\begin{eqnarray}
u_{m}(v_{-1}w)+F_{r+s+k-m+1}=v_{-1}(u_{m}w)+(u_{m}v)_{-1}w+F_{r+s+k-m+1}.
\end{eqnarray}
This proves $Y_{-}(u,x)\in x^{-1}({\rm Der} (\gr_{\E}(V)))[x^{-1}]$.
Therefore, $\gr_{\F}(V)$ is a vertex Poisson algebra.
\end{proof}

In the following, for each vertex algebra we construct 
a canonical decreasing sequence $\E=\{E_{n}\}_{n\ge 0}$ which satisfies
all the conditions assumed in Proposition \ref{pabstract}.
 
\bd{dsequence}
{\em Let $V$ be a vertex algebra and let $W$ be a $V$-module.
Define a sequence $\E_{W}=\{E_{n}(W)\}_{n\in \Z}$ of subspaces of $W$,
where for $n\in \Z$, $E_{n}(W)$ is linearly spanned by the vectors
\begin{eqnarray}\label{een-span}
u^{(1)}_{-1-k_{1}}\cdots u^{(r)}_{-1-k_{r}}w
\end{eqnarray}
for $r\ge 1,\; u^{(1)},\dots, u^{(r)}\in V,\;w\in W,\; k_{1},\dots,k_{r}\ge 0$ with
$k_{1}+\cdots +k_{r}\ge n$.}
\ed

Our main task is to establish the properties 
(\ref{edfil-condition}) and (\ref{evpacondition}) for the sequence $\E$.
The following are some immediate consequences:

\bl{leasy}
For any $V$-module $W$ we have
\begin{eqnarray}
& &E_{n}(W)\supset E_{n+1}(W)\;\;\;\mbox{ for any }n\in \Z,
\label{edecreasing}\\
& &E_{n}(W)=W\;\;\;\mbox{ for any }n\le 0,\label{enegative}\\
& &u_{-1-k}E_{n}(W)\subset E_{n+k}(W)\;\;\;\mbox{ for }u\in V,\; k\ge 0,\;
n\in \Z.\label{einclusion}
\end{eqnarray}
\el 

The following gives a stronger spanning property for $E_{n}(W)$:

\bl{lgeneral-ele}
Let $W$ be a $V$-module. For any $n\ge 1$,  we have
\begin{eqnarray}\label{e-en+1}
E_{n}(W)={\rm span} \{ u_{-1-i}w\;|\; u\in V,\; i\ge 1,\; w\in E_{n-i}(W)\}.
\end{eqnarray}
Furthermore, for $n\ge 1$, $E_{n}(W)$ is linearly
spanned by the vectors
\begin{eqnarray}\label{estronger-span}
u^{(1)}_{-k_{1}-1}\cdots u^{(r)}_{-k_{r}-1}w
\end{eqnarray}
for $r\ge 1,\; u^{(1)},\dots,u^{(r)}\in V,\;w\in W,\;
k_{1},\dots,k_{r} \ge 1$ with $k_{1}+\cdots +k_{r}\ge n$.
\el

\begin{proof} Notice that (\ref{estronger-span}) follows from 
(\ref{e-en+1}) and induction. 
Denote by $E_{n}'(W)$ the space on the right-hand side of (\ref{e-en+1}).
To prove (\ref{e-en+1}), we need to prove that each spanning vector of $E_{n}(W)$ 
in (\ref{een-span}) lies in $E_{n}'(W)$.
Now we use induction on $r$.
If $r=1$, we have $k_{1}\ge n\ge 1$ and $w\in W=E_{n-k_{1}}(W)$, 
so that $u^{(1)}_{-1-k_{1}}w\in E_{n}'(W)$.
Assume $r\ge 2$. If $k_{1}\ge 1$, we have
$$u^{(1)}_{-1-k_{1}}u^{(2)}_{-1-k_{2}}\cdots u^{(r)}_{-1-k_{r}}w\in E_{n}'(W)$$
because $u^{(2)}_{-1-k_{2}}\cdots u^{(r)}_{-1-k_{r}}w\in E_{n-k_{1}}(W)$
with $k_{2}+\cdots +k_{r}\ge n-k_{1}$.
If $k_{1}=0$, we have $k_{2}+\cdots +k_{r}\ge n$, so that
$u^{(2)}_{-1-k_{2}}\cdots u^{(r)}_{-1-k_{r}}w\in E_{n}(W)$.
By the inductive hypothesis, we have
$u^{(2)}_{-1-k_{2}}\cdots u^{(r)}_{-1-k_{r}}w\in E_{n}'(W)$.
Furthermore, for any $b\in V,\; k\ge 1,\; w'\in E_{n-k}(W)$, we have
$$u^{(1)}_{-1}b_{-1-k}w'=b_{-1-k}u^{(1)}_{-1}w'+\sum_{i\ge 0}{-1\choose i}(u^{(1)}_{i}b)_{-2-k-i}w'.$$
{}From definition we have $u^{(1)}_{-1}w'\in  E_{n-k}(W)$, so that
$b_{-1-k}u^{(1)}_{-1}w'\in  E_{n}'(W)$.
On the other hand, for $i\ge 0$, we have
$w'\in E_{n-k}(W)\subset E_{n-k-i-1}(W)$, so that
$$(u^{(1)}_{i}b)_{-2-k-i}w'\in E_{n}'(W).$$
Therefore, $u^{(1)}_{-1}b_{-1-k}w'\in E_{n}'(W).$
This proves that 
$u^{(1)}_{-1-k_{1}}u^{(2)}_{-1-k_{2}}\cdots u^{(r)}_{-1-k_{r}}w\in E_{n}'(W)$,
completing the induction.
\end{proof}

We have the following special case of
(\ref{edfil-condition}) and (\ref{evpacondition}) for $\E$:

\bl{lbaasicproperty}
Let $W$ be any $V$-module. For $a\in V,\; m,n\in \Z$, we have
\begin{eqnarray}
a_{m}E_{n}(W)\subset E_{n-m-1}(W).\label{enge-0}
\end{eqnarray}
Furthermore,
\begin{eqnarray}
a_{m}E_{n}(W)\subset E_{n-m}(W)\;\;\;\mbox{ for }m\ge 0.\label{enonneg}
\end{eqnarray}
\el

\begin{proof} By (\ref{einclusion}), (\ref{enge-0}) holds for $m\le -1$.
Assume $m\ge 0$. Since $E_{n-m}(W)\subset E_{n-m-1}(W)$
it suffices to prove (\ref{enonneg}).
We now prove the assertion by induction on $n$.
If $n\le 0$,
we have $E_{n-m}(W)=W$ (because $n-m\le 0$), so that
$a_{m}E_{n}(W)\subset W=E_{n-m}(W)$. 
Assume $n\ge 1$. From (\ref{e-en+1}), $E_{n}(W)$ is spanned by
the vectors $u_{-1-k}w$ for $u\in V, \; k\ge 1,\; w\in E_{n-k}(W)$.
Let $u\in V, \; k\ge 1,\; w\in E_{n-k}(W)$. In view of Borcherds' commutator
formula we have
$$a_{m}u_{-1-k}w
=u_{-1-k}a_{m}w+\sum_{i\ge 0}\binom{m}{i}(a_{i}u)_{m-k-i-1}w.$$
Since $w\in E_{n-k}(W)$ with $n-k<n$, from inductive hypothesis we have
\begin{eqnarray*}
& &a_{m}w\in E_{n-k-m}(W),\\
& &(a_{i}u)_{m-k-i-1}w\in E_{n-m+i}(W)\subset E_{n-m}(W)\;\;\;\mbox{ for }i\ge 0.
\end{eqnarray*}
Furthermore, using inductive hypothesis and Lemma \ref{leasy} we have
$$u_{-1-k}a_{m}w\in u_{-1-k}E_{n-k-m}(W)\subset E_{n-m}(W).$$
Therefore, $a_{m}u_{-1-k}w\in E_{n-m}(W)$. This proves $a_{m}E_{n}(W)\subset E_{n-m}(W)$,
completing the induction and the whole proof.
\end{proof}

Now we have the following general case:

\bp{pbasic2} 
Let $W$ be a $V$-module and let $u\in E_{r}(V),\; w\in E_{s}(W)$ 
with $r,s\in \Z$. Then
\begin{eqnarray}\label{eunw-general}
u_{n}w\in E_{r+s-n-1}(W)\;\;\;\mbox{ for }n\in \Z.
\end{eqnarray}
Furthermore, we have
\begin{eqnarray}\label{eunw-nonnegative}
u_{n}w\in E_{r+s-n}(W)\;\;\;\ \ \ \mbox{ for }n\ge 0.
\end{eqnarray}
\ep

\begin{proof} We are going to use induction on $r$. 
By Lemma \ref{lbaasicproperty}, we have
$u_{n}w\in E_{s-n-1}(W)$. If $r\le 0$, we have
$r+s-n-1\le s-n-1$, so that $u_{n}w\in E_{s-n-1}(W)\subset E_{r+s-n-1}(W)$.
Assume $r\ge 0$ and $u\in E_{r+1}(V)$. In view of (\ref{e-en+1})
it suffices to consider $u=a_{-2-i}b$ for some $a\in V,\; 0\le i\le r,\; b\in E_{r-i}(V)$.
By the iterate formula (\ref{eiterate}) we have
\begin{eqnarray}
(a_{-2-i}b)_{n}w
=\sum_{j\ge 0}(-1)^{j}\binom{-2-i}{j}
\left(a_{-2-i-j}b_{n+j}w-(-1)^{i}b_{n-2-i-j}a_{j}w\right).
\end{eqnarray}
If $n\ge 0$, using the inductive hypothesis (with $b\in E_{r-i}(V)$) and 
Lemma \ref{lbaasicproperty} we have
\begin{eqnarray*}
& &a_{-2-i-j}b_{n+j}w\in a_{-2-i-j}E_{r-i+s-n-j}(W)\subset E_{r+1+s-n}(W),\\
& &b_{n-2-i-j}a_{j}w\in b_{n-2-i-j}E_{s-j}(W)\subset E_{r+1+s-n}(W),
\end{eqnarray*}
{}from which we have that $(a_{-2-i}b)_{n}w\in E_{r+1+s-n}(W)$.
If $n\le -1$, we have
\begin{eqnarray*}
& &a_{-2-i-j}b_{n+j}w\in a_{-2-i-j}E_{r-i+s-n-j-1}(W)\subset E_{r+s-n}(W),\\
& &b_{n-2-i-j}a_{j}w\in b_{n-2-i-j}E_{s-j}(W)
\subset E_{r+1+s-n}(W)\subset E_{r+s-n}(W),
\end{eqnarray*}
so that $(a_{-2-i}b)_{n}w\in E_{r+s-n}(W)=E_{(r+1+s)-n-1}(W)$.
This concludes the proof.
\end{proof}

Combining Propositions \ref{pabstract} and \ref{pbasic2} we immediately have:

\bt{tmain1}
Let $V$ be any vertex algebra and let $\E=\{E_{n}(V)\}$ 
be the decreasing sequence 
defined in Definition \ref{dsequence} for $V$. Set
\begin{eqnarray}
\gr_{\E}(V)=\coprod_{n\ge 0} E_{n}/E_{n+1}.
\end{eqnarray}
Then $\gr_{\E}(V)$ equipped with the multiplication defined by
\begin{eqnarray}
(a+E_{r+1})(b+E_{s+1})=a_{-1}b+E_{r+s+1}
\end{eqnarray}
is a commutative and associative $\N$-graded algebra 
with ${\bf 1}+E_{1}\in E_{0}/E_{1}$
as identity and with a derivation $\partial$ defined by
\begin{eqnarray}
\partial (u+E_{n+1})=\D (u)+E_{n+2}\;\;\;\mbox{ for }u\in E_{n},\; n\in \N. 
\end{eqnarray}
Furthermore, $\gr_{\E}(V)$ is a vertex Poisson algebra where
\begin{eqnarray}
Y_{-}(a+E_{r+1},x)(b+E_{s+1})=\sum_{n\ge 0}(u_{n}v+E_{r+s-n+1})x^{-n-1}
\end{eqnarray}
for $a\in E_{r},\; b\in E_{s}$ with $r,s\in \N$.
\et

\bp{pmodule-vpa}
Let $W$ be any $V$-module and $\E_{W}$ 
the decreasing sequence defined in Definition \ref{dsequence} for $W$.
Then the associated graded vector space 
$\gr_{\E}(W)=\coprod_{n\ge 0}E_{n}(W)/E_{n+1}(W)$ 
is naturally a module for the vertex Poisson algebra $\gr_{\E}(V)$
with
\begin{eqnarray}.
& &(v+E_{r+1}(V))\cdot (w+E_{s+1}(W))=v_{-1}w+E_{r+s+1}(W),
\label{ecaa-action}\\
& &Y_{-}^{W}(v+E_{r+1})(w+E_{s+1}(W))
=\sum_{n\ge 0}(v_{n}w+E_{r+s-n+1})x^{-n-1}
\label{evertex-action}
\end{eqnarray}
for $v\in E_{r}(V),\; w\in E_{s}(W)$.
\ep

\begin{proof} With the properties (\ref{eunw-general}) 
and (\ref{eunw-nonnegative})
the actions given by (\ref{ecaa-action}) and (\ref{evertex-action}) 
are well defined.
Clearly, ${\bf 1}+E_{1}$ acts on $\gr_{\E}(W)$ as identity and
we have
$$(u+E_{r+1}(V))\cdot \left( (v+E_{s+1}(V))\cdot (w+E_{k+1}(W))\right)
=u_{-1}v_{-1}w+E_{r+s+k+1}(W).$$
By the iterate formula (\ref{eiterate}) we have
$$(u_{-1}v)_{-1}w
=\sum_{i\ge 0}\left(u_{-1-i}v_{-1+i}w+v_{-2-i}u_{i}w\right),$$
where for $i\ge 1$, using (\ref{eunw-nonnegative}) and 
(\ref{eunw-general}) we have
$$u_{-1-i}v_{-1+i}w\in u_{-i-1}E_{s+k+1-i}(W)\subset E_{r+s+k+1}(W)$$
and for $i\ge 0$, similarly we have
$$v_{-2-i}u_{i}w\in v_{-2-i}E_{s+k-i}(W)\subset E_{r+s+k+1}(W).$$
Thus
$$(u_{-1}v)_{-1}w\in u_{-1}v_{-1}w+E_{r+s+k+1}(W).$$
This proves that $\gr_{\E}(W)$ is a module for $\gr_{\E}(V)$ 
as an associative algebra.
It is straightforward to check that 
it is a module for the vertex Lie algebra. 
Other properties are clear from the proof of Proposition \ref{pabstract}.
\end{proof}

Notice that so far we have not excluded the possibility that the
associated sequence $\E_{V}$ is {\em trivial} in the sense that $E_{n}(V)=V$
for all $n\ge 0$.  Indeed, as we shall see in the next section, for
some vertex algebras the associated sequence $\E$ is trivial. 

Nevertheless, we have:

\bl{lNgva}
Let $V=\coprod_{n\ge 0}V_{(n)}$ be an $\N$-graded vertex algebra 
and $\E=\{E_{n}\}$ be the decreasing sequence 
defined in Definition \ref{dsequence} for $V$. Then 
\begin{eqnarray}\label{engva}
E_{n}(V)\subset \coprod_{m\ge n}V_{(m)}\;\;\;\mbox{ for }n\ge 0.
\end{eqnarray}
Furthermore, the associated decreasing sequence $\E=\{ E_{n}\}$ 
for $V$ is a filtration, i.e.,
\begin{eqnarray}\label{enva-fil}
\cap_{n\ge 0}E_{n}(V)=0.
\end{eqnarray}
\el

\begin{proof} By definition we have $E_{0}=V=\coprod_{n\ge 0}V_{(n)}$.
For $n\ge 1$, recall that $E_{n}$ is
linearly spanned by the vectors
\begin{eqnarray*}
u^{(1)}_{-1-k_{1}}\cdots u^{(r)}_{-1-k_{r}}v
\end{eqnarray*}
for $r\ge 1,\; u^{(1)},\dots, u^{(r)},v\in V,\; k_{1},\dots,k_{r}\ge 1$ with
$k_{1}+\cdots +k_{r}\ge n$. 
If the vectors $u^{(1)},\dots, u^{(r)},v$ are homogeneous, we have
$$\wt\left(u^{(1)}_{-1-k_{1}}\cdots u^{(r)}_{-1-k_{r}}v\right)
=\wt u^{(1)}+k_{1}+\cdots+ \wt u^{(r)}+k_{r}+\wt v
\ge k_{1}+\cdots +k_{r}\ge n.$$
This proves (\ref{engva}) for $n\ge 1$. Clearly,
each subspace $E_{n}$ of $V$ is graded. From (\ref{engva})
we immediately have (\ref{enva-fil}).
\end{proof}

In the next section we shall generalize Lemma \ref{lNgva} from
an $\N$-graded vertex algebra to a lower truncated $\Z$-graded vertex algebra
by using a relation between the decreasing sequence $\E$ and
a sequence introduced by Zhu.

\section{The relation between the sequences $\E$ and
${\mathcal{C}}$}

In this section we first recall the sequence
 ${\mathcal{C}}$ introduced by Zhu and we then give a
relation between the two decreasing sequences $\E$ and
${\mathcal{C}}$.  We show
that if $V$ is a lower truncated $\Z$-graded vertex algebra, then both
sequences are decreasing filtrations of $V$. 

The following definition is (essentially) due to Zhu ([Z1,2]):

\bd{dcspaces}
{\em Let $V$ be a vertex algebra and $W$ a $V$-module.
For any $n\ge 2$ we define $C_{n}(W)$ to be
the subspace of $W$, linearly spanned by 
the vectors $v_{-n}w$ for $v\in V,\; w\in W$.
A $V$-module $W$ is said to be 
{\em $C_{n}$-cofinite} if $W/C_{n}(W)$ is finite-dimensional.
In particular, if $V/C_{n}(V)$ is finite-dimensional,
we say that the vertex algebra $V$ is {\em $C_{n}$-cofinite}.}
\ed 

The following are easy consequences:

\bl{lcn-basic}
Let $V$ be any vertex algebra, let $W$ be a $V$-module 
 and let $n\ge 2$. Then
\begin{eqnarray}
& &C_{m}(W)\subset C_{n}(W)\;\;\;\mbox{ for }m\ge n,\label{ecn-inclusion}\\
& &u_{-k}C_{n}(W)\subset C_{n}(W)
\;\;\;\mbox{ for }u\in V,\; k\ge 0,\label{e3.2}\\
& &u_{-n}v_{-k}w\equiv v_{-k}u_{-n}w\mod C_{n+k}(W)
\;\;\;\mbox{ for }u,v\in V,\;w\in W,\; k\ge 0.
\label{ecomm-uv}
\end{eqnarray}
\el

\begin{proof} For $v\in V,\; r\ge 2$ we have 
$v_{-r-1}=\frac{1}{r}(\D v)_{-r}$. 
{}From this we immediately have $C_{r+1}(W)\subset C_{r}(W)$ for $r\ge
2$, which implies (\ref{ecn-inclusion}).  Let $u,v\in V,\;w\in W,\;
k\ge 0$. Using the commutator formula (\ref{ecommutator}) and
(\ref{ecn-inclusion}) we have
$$u_{-k}v_{-n}w
=v_{-n}u_{-k}w+\sum_{i\ge 0}\binom{-k}{i}(u_{i}v)_{-k-n-i}w\in C_{n}(W),$$
proving that $u_{-k}C_{n}(W)\subset C_{n}(W)$. We also have
$$u_{-n}v_{-k}w-v_{-k}u_{-n}w
=\sum_{i\ge 0}\binom{-n}{i}(u_{i}v)_{-n-k-i}w\in C_{n+k}(W),$$
proving (\ref{ecomm-uv}). 
\end{proof}

We also have the following more technical results:

\bl{lmorecn}
Let $V$ be any vertex algebra, let $W$ be a $V$-module and let $k\ge 2$. Then
\begin{eqnarray}
u_{-k}C_{k}(W)\subset C_{k+1}(W)\;\;\;\mbox{ for }u\in V.
\end{eqnarray}
\el

\begin{proof} For $u,v\in V,\;w\in W$, in view of the iterate formula (\ref{eiterate}) we have
\begin{eqnarray}\label{eu-1v-2k+1}
(u_{-1}v)_{-2k+1}w
=\sum_{i\ge 0}\left(u_{-1-i}v_{-2k+1+i}w+v_{-2k-i}u_{i}w\right).
\end{eqnarray}
Now we examine each term in (\ref{eu-1v-2k+1}).
Notice that $(u_{-1}v)_{-2k+1}w\in C_{k+1}(W)$ as $-2k+1\le -k-1$ and that
$v_{-2k-i}u_{i}w\in C_{k+1}(W)$ for $i\ge 0$ as $-2k-i\le -k-1$.
If $i\ge k$, we have $-1-i\le -k-1$, so that 
$u_{-1-i}v_{-2k+1+i}w\in C_{k+1}(W)$.
For $0\le i\le k-2$, we have $-2k+1+i\le -k-1$, so that
$v_{-2k+1+i}w\in C_{k+1}(W)$.
Then by Lemma \ref{lcn-basic} we have $u_{-1-i}v_{-2k+1+i}w\in C_{k+1}(W)$
for $0\le i\le k-2$. Therefore, the only remaining term 
$u_{-k}v_{-k}w$ in (\ref{eu-1v-2k+1}) must also lie
in $C_{k+1}(W)$. This proves $u_{-k}C_{k}(W)\subset C_{k+1}(W)$.
\end{proof}

\bp{pcn}
Let $V$ be any vertex algebra, let $W$ be a $V$-module and 
let $n$ be any nonnegative integer. Then
\begin{eqnarray}
u^{(1)}_{-k_{1}}\cdots u^{(r)}_{-k_{r}}w\in C_{n+2}(W)
\end{eqnarray}
for $r\ge 2^{n},\; u^{(1)},\dots,u^{(r)}\in V,\;w\in W,
\;k_{1},\dots,k_{r}\ge 2$.
\ep

\begin{proof} Since $u_{-i}C_{n+2}(W)\subset C_{n+2}(W)$ for $u\in V,\; i\ge 0$
(by Lemma \ref{lcn-basic}), it suffices to prove the assertion for $r=2^{n}$.
Also, since $u_{-k}=\frac{1}{(k-1)!} (\D^{k-2} u)_{-2}$ for $u\in V,\; k\ge 2$,
it suffices to prove the assertion for $k_{1}=\cdots =k_{r}=2$.
We are going to use induction on $n$. If $n=0$,
by definition we have $v_{-2}w\in C_{2}(W)$ for $v\in V,\; w\in W$.
Assume the assertion holds for $n=p$, some nonnegative integer.
Assume that $r=2^{p+1}$ and set $s=2^{p}$. 
Let $u^{(1)},\dots, u^{(r)}\in V,\; w\in W$.
By inductive hypothesis we have
\begin{eqnarray*}
u^{(s+1)}_{-2}\cdots u^{(r)}_{-2}w\in C_{p+2}(W),
\end{eqnarray*}
so that
\begin{eqnarray}\label{e3.8}
u^{(1)}_{-2}\cdots u^{(r)}_{-2}w\in 
u^{(1)}_{-2}\cdots u^{(s)}_{-2} C_{p+2}(W).
\end{eqnarray}
Consider a typical spanning vector $a_{-p-2}w'$ of $C_{p+2}(W)$
for $a\in V,\; w'\in W$. Using (\ref{ecomm-uv}) and (\ref{e3.2}) we have
\begin{eqnarray}\label{e3.9}
u^{(1)}_{-2}\cdots u^{(s)}_{-2} a_{-p-2}w'
\equiv a_{-p-2}u^{(1)}_{-2}\cdots u^{(s)}_{-2}w'\mod C_{p+4}(W).
\end{eqnarray}
Furthermore, by inductive hypothesis, we have
$$u^{(1)}_{-2}\cdots u^{(s)}_{-2}w'\in C_{p+2}(W),$$
which together with Lemma \ref{lmorecn} gives
\begin{eqnarray}
a_{-p-2}u^{(1)}_{-2}\cdots u^{(s)}_{-2}w'\in a_{-p-2}C_{p+2}(W)
\subset C_{p+3}(W).
\end{eqnarray}
Thus by (\ref{e3.9}) we have
$$u^{(1)}_{-2}\cdots u^{(s)}_{-2} a_{-p-2}w'\in C_{p+3}(W),$$
proving that
\begin{eqnarray}
u^{(1)}_{-2}\cdots u^{(s)}_{-2} C_{p+2}(W)\subset C_{p+3}(W).
\end{eqnarray}
Therefore, by (\ref{e3.8}) we have
$$u^{(1)}_{-2}\cdots u^{(r)}_{-2}w\in C_{p+3}(W).$$
This finishes the induction steps and completes the proof.
\end{proof}

The relation between the two decreasing sequences $\{E_{n}(W)\}$ and 
$\{C_{n}(W)\}$ is described as follows:

\bt{tZlower}
Let $W$ be any module for vertex algebra $V$
and let $\E_{W}=\{E_{n}(W)\}$ be the associated decreasing sequence. Then
for any $n\ge 2$,
\begin{eqnarray}
& &C_{n}(W)\subset E_{n-1}(W),\label{ecn-in-en-1}\\
& &E_{m}(W)\subset C_{n}(W)\;\;\;\mbox{ whenever }m\ge (n-2)2^{n-2}.
\label{eem-in-cn}
\end{eqnarray}
Furthermore,
\begin{eqnarray}\label{esame-inter}
\cap_{n\ge 0}E_{n}(W)=\cap_{n\ge 0}C_{n+2}(W).
\end{eqnarray}
\et

\begin{proof}  {}From the definitions of $C_{n}(W)$ and $E_{n-1}(W)$
we immediately have $C_{n}(W)\subset E_{n-1}(W)$.
Consider a generic spanning element of $E_{m}(W)$:
$$X=u^{(1)}_{-1-k_{1}}\cdots u^{(r)}_{-1-k_{r}}w$$
where $r\ge 1,\; u^{(1)},\dots,u^{(r)}\in V,\;w\in W,
\; k_{1},\dots,k_{r}\ge 1$ 
with $k_{1}+\cdots +k_{r}\ge m$.
If $k_{i}\ge n-1$ for some $i$, by
(\ref{ecn-inclusion}) we have
$u^{(i)}_{-1-k_{i}}W\subset C_{-1-k_{i}}(W)\subset C_{n}(W)$
and then by (\ref{e3.2}) we have $X\in C_{n}(W)$.
If $r\ge 2^{n-2}$, by Proposition \ref{pcn} $X\in C_{n}(W)$.
Since $k_{1}+\cdots +k_{r}\ge m\ge (n-2)2^{n-2}$,
we have either $k_{i}\ge n-1$ for some $i$ or $r\ge 2^{n-2}$.
Therefore, $X\in C_{n}(W)$ whenever $m\ge (n-2)2^{n-2}$. 
This proves (\ref{eem-in-cn}).
Combining (\ref{eem-in-cn}) and (\ref{ecn-in-en-1}) 
we have (\ref{esame-inter}).
\end{proof}

\bc{ce1e2}
For any vertex algebra $V$ and any $V$-module $W$, we have
\begin{eqnarray}
E_{1}(W)=C_{2}(W),\ \ \  E_{2}(W)=C_{3}(W).
\end{eqnarray}
\ec

\begin{proof} By (\ref{ecn-in-en-1}) we have $C_{2}(W)\subset E_{1}(W)$ and 
$C_{3}(W)\subset E_{2}(W)$. On the other hand,
by (\ref{eem-in-cn}) with $m=1, n=2$ we have
$E_{1}(W)\subset C_{2}(W)$ and by (\ref{eem-in-cn}) with $m=2, n=3$
we have $E_{2}(W)\subset C_{3}(W)$.
\end{proof}

Recall the following result of Zhu [Z1,2]:

\bp{pzh}
Let $V$ be any vertex algebra. Then $V/C_{2}(V)$
is a Poisson algebra with
\begin{eqnarray}
\bar{u}\cdot \bar{v}=\overline{u_{-1}v},\;\;\;\;
[\bar{u},\bar{v}]=\overline{u_{0}v}\;\;\;\;\mbox{ for }u,v\in V,
\end{eqnarray}
where $\bar{u}=u+C_{2}(V)$, and with ${\bf 1}+C_{2}(V)$ 
as the identity element.
\ep

It is clear that the degree zero subspace $E_{0}/E_{1}$ of $\gr_{\E}(V)$ is
a Poisson algebra where
$$(u+E_{1}) (v+E_{1})=u_{-1}v+E_{1},\;\;\;\;
[u+E_{1},v+E_{1}]=u_{0}v+E_{1}$$
for $u,v\in V$.
With $E_{0}(V)=V$ and $E_{1}(V)=C_{2}(V)$, we see that
this Poisson algebra is nothing but the
Zhu's Poisson algebra $V/C_{2}(V)$.

Thus we have:

\bp{pzhu-degree0}
Let $V$ be any vertex algebra.
The degree zero subspace 
$\gr_{\E}(V)_{(0)}=E_{0}(V)/E_{1}(V)$ of
the $\N$-graded vertex Poisson algebra $\gr_{\E}(V)$ is 
naturally a Poisson algebra which coincides with the
Zhu's Poisson algebra $V/C_{2}(V)=E_{0}/E_{1}$.
\ep

The following result generalizes the result of Lemma \ref{lNgva}:

\bp{pz-gradedva}
Let $V=\coprod_{n\ge t}V_{(n)}$ be a lower truncated $\Z$-graded 
vertex algebra. Then
\begin{eqnarray}\label{ecn+2-in}
C_{n}(V)\subset \coprod_{k\ge 2t+n-1}V_{(k)}
\end{eqnarray}
for $n\ge 2$.
Furthermore,
\begin{eqnarray}
& &E_{m}(V)\subset \coprod_{k\ge 2t+n-1}V_{(k)}
\;\;\;\mbox{ whenever }m\ge (n-2)2^{n-2},
\label{ekey-need}\\
& &\cap_{n\ge 0}E_{n}(V)=\cap_{n\ge 2} C_{n}(V)=0.\label{e-cn-int=0}
\end{eqnarray}
\ep

\begin{proof} For homogeneous vectors $u,v\in V$ and for any $n\ge 2$  we have
$$\wt (u_{-n}v)=\wt u+\wt v+n-1\ge 2t+n-1.$$
In view of this we have
$$C_{n}(V)\subset \coprod_{k\ge 2t+n-1}V_{(k)}$$
for $n\ge 2$. This proves (\ref{ecn+2-in}), from which we immediately have 
that $\cap_{n\ge 0} C_{n+2}(V)=0$.
Using Theorem \ref{tZlower}, we obtain
(\ref{ekey-need}) and (\ref{e-cn-int=0}).
\end{proof}

For the rest of this section, we consider vertex algebras whose associated
decreasing sequence $\E$ is trivial.

First we have:

\bl{ltrivial}
Let $V$ be a vertex algebra and let $W$ be a $V$-module. If $W=C_{2}(W)$, then
\begin{eqnarray}
E_{n}(W)=C_{n+2}(W)=W\;\;\;\mbox{  for all }n\ge 0.
\end{eqnarray}
\el

\begin{proof}
Since $W=C_{2}(W)$, we have $E_{1}(W)=C_{2}(W)=W$.
Assume that $E_{k}(W)=W$ for some $k\ge 1$.
Then
$$v_{-2}W=v_{-2}E_{k}(W)\subset E_{k+1}(W)\;\;\;\mbox{ for }v\in V.$$
{}From this we have $W=C_{2}(W)\subset E_{k+1}(W)$, proving
$E_{k+1}(W)=W$. By induction, we have $E_{n}(W)=W$ for all $n\ge 0$.
In view of Theorem \ref{tZlower} we have $C_{n}(W)=W$ for all $n\ge 2$.
\end{proof}

Suppose that $V$ is a vertex algebra such that $C_{2}(V)=V$.
By Lemma \ref{ltrivial} we have $C_{n+2}(V)=V$ for $n\ge 0$, so that
$V=\cap_{n\ge 0}C_{n+2}(V)$. Furthermore,
if there exists a lower truncated $\Z$-grading $V=\coprod_{n\in \Z}V_{(n)}$
with which $V$ becomes an $\Z$-graded vertex algebra, by (\ref{e-cn-int=0})
(Proposition \ref{pz-gradedva}) we have $\cap_{n\ge 0}C_{n+2}(V)=0$, so that
$V=\cap_{n\ge 0}C_{n+2}(V)=0$.
Therefore we have proved:

\bp{pnonexistence}
Let $V$ be a nonzero vertex algebra such that $C_{2}(V)=V$. 
Then there does not exist a lower truncated $\Z$-grading 
$V=\coprod_{n\in \Z}V_{(n)}$
with which $V$ becomes an $\Z$-graded vertex algebra.
\ep

{}From [B] and [FLM], associated to any
nondegenerate even lattice $L$ of finite rank, we have a vertex algebra
$V_{L}$. Furthermore, $V_{L}$ is a vertex operator algebra if and only
if $L$ is positive-definite in the sense that $\<\alpha,\alpha\>>0$ 
for $0\ne \alpha\in L$. 
In this case, $V_{L}$ is $\N$-graded by
$L(0)$-weight (with $1$-dimensional weight-zero subspace), 
so that Lemma \ref{lNgva} (and Proposition \ref{pz-gradedva}) 
applies to $V_{L}$. On the other hand, we have:

\bp{platticevoa} 
Let $L$ be a finite rank nondegenerate even lattice that is not positive-definite
and let $V_{L}$ be the associated vertex algebra.  Then
$C_{n+2}(V_{L})=E_{n}(V_{L})=V_{L}$ for $n\ge 0$.
Furthermore, there does not exist a lower truncated $\Z$-grading on $V_{L}$ with which
$V_{L}$ becomes a lower truncated $\Z$-graded vertex algebra.
\ep

\begin{proof} First we show that there exists $\alpha\in L$ such that $\<\alpha,\alpha\><0$.
Since $L$ is not positive-definite, there exists $0\ne \beta\in L$ such that
$\<\beta,\beta\>\le 0$. If $\<\beta,\beta\>\ne 0$, that is,
$\<\beta,\beta\><0$, then we can simply take $\alpha=\beta$.
Suppose $\<\beta,\beta\>= 0$. Since $L$ is nondegenerate, there exists
$\gamma\in L$ such that $\<\gamma,\beta\>\ne 0$.
For $m\in \Z$, we have
$$\<\gamma+m\beta,\gamma+m\beta\>=\<\gamma,\gamma\>+2m\<\gamma,\beta\>.$$
We see that $\<\gamma+m\beta,\gamma+m\beta\><0$ for some $m$. Then
we can take $\alpha=\gamma+m\beta$ with the desired property.

Let $\alpha\in L$ be such that $\<\alpha,\alpha\><0$ and set
$\<\alpha,\alpha\>=-2k$ with $k\ge 1$.
Using the explicit expression of the vertex operators in [FLM],
we have $(e^{\alpha})_{-2k-1} e^{-\alpha}={\bf 1}$, so that 
${\bf 1}\in C_{2k+1}(V_{L})\subset C_{2}(V_{L})$.
Then $v=v_{-1}{\bf 1}\in C_{2}(V_{L})$ for $v\in V_{L}$. 
Thus $C_{2}(V_{L})=V_{L}$. 
By Lemma \ref{ltrivial} we have 
$C_{n+2}(V_{L})=E_{n}(V_{L})=V_{L}$ for $n\ge 0$.
The last assertion follows immediately from Proposition \ref{pnonexistence}.
\end{proof}

\section{Generating subspaces of vertex algebras}
In this section we shall use the differential
algebra structure on $\gr_{\E}(V)$ to study certain kinds of
generating subspaces of lower truncated $\Z$-graded vertex algebras.

First we prove the following results for classical algebras:

\bl{ldiff}
Let $(A,\partial)$ be an $\N$-graded (unital) differential algebra
such that $(\partial A)A=A_{+}$, where 
\begin{eqnarray}
A_{+}=\coprod_{n\ge 1}A_{(n)}.
\end{eqnarray}
Let $S$ be a generating subspace of $A_{(0)}$ as an algebra. Then
$A$ is linearly spanned by the vectors
\begin{eqnarray}
\partial^{n_{1}}(a_{1})\cdots \partial^{n_{r}}(a_{r})
\end{eqnarray}
for $r\ge 0,\; n_{1}\ge n_{2}\ge \cdots\ge n_{r}\ge 0,\; a_{1},\dots,a_{r}\in
S$, or equivalently, $S$ generates $A$ as a differential algebra. 
In particular, $A_{(0)}$ generates $A$ as a differential algebra.
Furthermore, $A$ is linearly spanned by the vectors
\begin{eqnarray}
\partial^{n_{1}}(a_{1})\cdots \partial^{n_{r}}(a_{r})
\end{eqnarray}
for $r\ge 1,\; n_{1}>n_{2}>\cdots> n_{r}\ge 0,\; a_{1},\dots,a_{r}\in
A_{(0)}$. 
\el

\begin{proof} First, we show that $A$ as a differential algebra is
generated by $A_{(0)}$.
Let $A'$ be the differential subalgebra of $A$, generated by
$A_{(0)}$. We are going to show (by induction) 
 that $\coprod_{n=0}^{k}A_{(n)}\subset A'$ for all $k\ge 0$.
{}From definition, we have $A_{(0)}\subset A'$.
Assume that $\coprod_{n=0}^{k}A_{(n)}\subset A'$ for some $k\ge 0$.
Consider the subspace $A_{(k+1)}$ of $A$. From our assumption,
we have
$$A_{(k+1)}\subset A_{+}=A\partial A,$$
so $A_{(k+1)}$ is linearly spanned by the vectors
$a\partial b$ for $a\in A_{(r)},\; b\in A_{(s)}$ with $r+s+1=k+1$.
For any $a\in A_{(r)},\; b\in A_{(s)}$ with $r+s+1=k+1$,
since $r,s\le k$ (with $r,s\ge 0$), by the inductive hypothesis,
we have $a,b\in A'$. Consequently, $a\partial b\in A'$.
Thus $A_{(k+1)}\subset A'$. This proves that
$\coprod_{n=0}^{k}A_{(n)}\subset A'$ for all $k\ge 0$.
Therefore, we have $A=A'$, proving that
$A$ as a differential algebra is generated by $A_{(0)}$.
It follows that if $S$ generates $A_{(0)}$ as an algebra,
then $S$ generates $A$ as a differential algebra.

For a positive integer $n$, let $A_{(n)}''$ be the subspace of $A_{(n)}$ 
spanned by the vectors
\begin{eqnarray}\label{e3.20}
\partial^{k_{1}}(a_{1})\cdots \partial^{k_{r}}(a_{r})b
\end{eqnarray}
for $r\ge 1,\; k_{1}>k_{2}>\cdots> k_{r}\ge 1,\; a_{1},\dots,a_{r},b\in
A_{(0)}$ with $k_{1}+\cdots +k_{r}=n$. We must prove 
$A_{(n)}=A_{(n)}''$ for all $n\ge 1$.

For a positive integer $n$, denote by $P_{n}$ 
the set of partitions of $n$. We now endow $P_{n}$ 
with the reverse order of the lexicographic order
on $P_{n}$. Set $P=\cup_{n\ge 1}P_{n}$.
For $\alpha\in P_{m},\;\beta\in P_{n}$, combining $\alpha$ and $\beta$
together we get 
a partition of $m+n$, which we denote by $\alpha*\beta$. Clearly, this
defines an abelian semigroup structure on $P$. Furthermore,
for $\alpha,\beta\in P_{n},\; \gamma\in P$, if $\alpha>\beta$, then
$\alpha*\gamma > \beta*\gamma$.
That is, the order is compatible with the multiplication.

For $\alpha\in P_{n}$, define $A_{(n)}^{\alpha}$ 
to be the linear span of the vectors
$$\partial^{k_{1}}(a_{1})\cdots \partial^{k_{r}}(a_{r})b$$
for $r\ge 1,\; k_{1}\ge k_{2}\ge \cdots \ge k_{r}\ge 1,\; a_{1},\dots,a_{r},b\in
A_{(0)}$ with $k_{1}+\cdots +k_{r}=n$ and
$(k_{1},\dots,k_{r})\le \alpha$. 
Since $A_{(0)}$ generates $A$ as a differential algebra, 
$\{ A_{(n)}^{\alpha}\}$ is a (finite) increasing filtration of $A_{(n)}$.

For $a,b\in A_{(0)}$ and $k\ge 1$, we have
\begin{eqnarray*}
\partial^{2k} (ab)=\sum_{i=0}^{2k}
\binom{2k}{i} \partial^{2k-i}(a)\partial^{i}(b),
\end{eqnarray*}
which can be rewritten as
\begin{eqnarray}\label{e3.24}
& &\binom{2k}{k}\partial^{k}(a)\partial^{k}(b)\nonumber\\
&=&\partial^{2k} (ab)-\partial^{2k} (a)b-\partial^{2k} (b)a-\sum_{i=1}^{k-1}
\binom{2k}{i} \left(\partial^{2k-i}(a)\partial^{i}(b)
+\partial^{i}(a)\partial^{2k-i}(b)\right).
\end{eqnarray}
We see that $(k,k)>(2k), (2k-i, i)$ for $1\le i\le k-1$.

Now consider a typical element of $A_{(n)}^{\alpha}$
$$X=\partial^{k_{1}}(a_{1})\cdots \partial^{k_{r}}(a_{r}) b$$
for $(k_{1}, \dots, k_{r})\in P_{n},\; a_{1},\dots,a_{r},b\in A_{(0)}$ 
with $(k_{1},\dots,k_{r})\le \alpha$. 
If all $k_{1}> k_{2}>\cdots >k_{r}$, then $X\in A_{(n)}''$.
Otherwise, using (\ref{e3.24}) we see that
$$X\in \sum_{\beta <\alpha}A_{(n)}^{\beta}.$$
Now it follows immediately from induction that $A_{(n)}=A_{(n)}''$.
\end{proof}

\bl{lcheck-assumption} 
Let $V$ be a vertex algebra and let $A=\gr_{\E}(V)$ be
the vertex Poisson algebra, obtained in Theorem \ref{tmain1}, which
is in particular an $\N$-graded (unital) differential algebra. 
Then $A_{+}=A\partial A$. Furthermore, for any $V$-module $W$,
the associated graded vector space 
$\gr_{\E}(W)$ is an $A$-module with
\begin{eqnarray}
(u+E_{m+1}(V))\cdot (w+E_{n+1}(W))=u_{-1}w+E_{m+n+1}(W)
\end{eqnarray}
for $u\in E_{m}(V),\; w\in E_{n}(W)$ with $m,n\in \N$, and
$\gr_{\E}(W)$ as an $A$-module is generated by
$E_{0}(W)/E_{1}(W)$, i.e., 
\begin{eqnarray}\label{emodule-top}
\gr_{\E}(W)=A(E_{0}(W)/E_{1}(W)).
\end{eqnarray}
\el

\begin{proof}
We have $A=\coprod_{n\in\N}A_{(n)}$, where $A_{(n)}=E_{n}/E_{n+1}$ for
$n\in \N$. For $n\ge 1$, from Lemma \ref{leasy},
$E_{n}$ is linearly spanned by the vectors
$u_{-2-i}v\in E_{n}$ where $u\in V,\; v\in E_{n-1-i}$ 
for $0\le i\le n-1$, and furthermore, we have
\begin{eqnarray*}
u_{-2-i}v+E_{n+1}&=&\frac{1}{(i+1)!}(\D^{i+1} u)_{-1}v+E_{n+1}\\
&=&\frac{1}{(i+1)!}(\D^{i+1} u+E_{i+2})(v+E_{n-i})\\
&=&\frac{1}{(i+1)!}\partial^{i+1} (u+E_{1}) (v+E_{n-i})\\
&\in&A\partial A,
\end{eqnarray*}
noticing that for any $r\in \Z$,
$\D E_{r}\subset E_{r+1}$ from the definition of $\D$ and Lemma \ref{pbasic2}.
This proves $E_{n}/E_{n+1}\subset A\partial A$ for $n\ge 1$, so that
$A_{+}\subset A\partial A$.
We also have that $A\partial A\subset AA_{+}\subset A_{+}$.
Therefore, $A\partial A=A_{+}$. 

For a $V$-module $W$, from Proposition \ref{pmodule-vpa}
$\gr_{\E}(W)$ is a module for $\gr_{\E}(V)$ as an algebra.
We must prove that 
$E_{n}(W)/E_{n+1}(W)\subset A(E_{0}(W)/E_{1}(W))$ for $n\ge 1$.
By Lemma \ref{leasy}, $E_{n}(W)$ is linearly spanned by the subspaces
$u_{-2-i}E_{n-1-i}(W)$ for $u\in V,\; 0\le i\le n-1$. 
For $w\in E_{n-1-i}(W)$, we have
\begin{eqnarray*}
u_{-2-i}w+E_{n+1}(W)=
\frac{1}{(i+1)!}\partial^{i+1} (u+E_{1}) (w+E_{n-i}(W)).
\end{eqnarray*}
Then it follows immediately from induction.
\end{proof}

Combining Lemmas \ref{ldiff} and \ref{lcheck-assumption} 
we immediately have:

\bc{cdiff}
Let $V$ be a vertex algebra and let $\gr_{\E}(V)$ be
the vertex Poisson algebra obtained in Theorem \ref{tmain1}. 
Then $\gr_{\E}(V)$ is linearly spanned by the vectors
\begin{eqnarray}
\partial^{k_{1}}(v^{(1)}+E_{1})\cdots
\partial^{k_{r}}(v^{(r)}+E_{1})
\end{eqnarray}
for $r\ge 1,\; v^{(i)}\in V,\; k_{1}>\cdots >k_{r}\ge 0$.
In particular, $\gr_{\E}(V)$ as a differential algebra
is generated by the subspace $E_{0}/E_{1}\;(=V/C_{2}(V))$.
\ec

The following result generalizes a theorem of [GN] (see also [NT]):

\bp{pfinite-prop}
Let $V$ be any vertex algebra. If $V$ is $C_{2}$-cofinite,
then $V$ is $E_{n}$-cofinite and $C_{n+2}$-cofinite for any $n\ge 0$.
\ep

\begin{proof} Since $\dim V/C_{2}<\infty$, it follows from 
Corollary \ref{cdiff} that for each $n\ge 0$, the degree $n$ subspace
$E_{n}/E_{n+1}$ of $\gr_{\E}(V)$ is finite dimensional.
Consequently, $\dim V/E_{n}=\dim E_{0}/E_{n}<\infty$ for all $n\ge 0$.
For any $n\ge 2$, by (\ref{ekey-need}) we have
$E_{m}\subset C_{n}$ for $m=(n-2)2^{n-2}$.
Then $\dim V/C_{n}\le \dim V/E_{m}<\infty$.
\end{proof}

Furthermore we have (cf. [Bu1,2]):

\bp{pmodule-cn}
Let $V$ be any vertex algebra and $W$ any $V$-module.
If $V$ and $W$ are $C_{2}$-cofinite, then $W$ is 
$C_{n}$-cofinite for all $n\ge 2$.
\ep

\begin{proof} In the proof of Proposition \ref{pfinite-prop}, 
we showed that
$\gr_{\E}(V)$ is an $\N$-graded differential algebra with
finite-dimensional homogeneous subspaces. Since $\dim
W/C_{2}(W)<\infty$, it follows from (\ref{emodule-top})
that all the homogeneous subspaces of $\gr_{\E}(W)$ are
finite-dimensional. The same argument of Proposition \ref{pfinite-prop}
shows that $W$ is $C_{n}$-cofinite for all $n\ge 2$.
\end{proof}

\br{rBuh}
{\em It has been proved in [Bu1] and [NT] that
if $V$ is a vertex operator algebra with nonnegative weights and with
$V_{(0)}=\C{\bf 1}$ and if $V$ is $C_{2}$-cofinite, then
any irreducible $V$-module $W$ is $C_{n}$-cofinite
for all $n\ge 2$.}
\er

The following result generalizes a theorem of [NG] (cf. [Bu1-2], [ABD]):

\bt{tmain-appl}
Let $V=\coprod_{n\ge t}V_{(n)}$ be any lower truncated $\Z$-graded 
vertex algebra
such as a vertex operator algebra in the sense of [FLM] and [FHL].
Then for any graded subspace $U$ of $V$,
$V=U+ C_{2}(V)$ if and only if $V$ is linearly spanned by the vectors
\begin{eqnarray}\label{emain-span}
u^{(1)}_{-n_{1}}\cdots u^{(r)}_{-n_{r}}{\bf 1}
\end{eqnarray}
for $r\ge 0,\; n_{1}>\cdots >n_{r}\ge 1,\; u^{(1)},\dots,u^{(r)}\in U$.
\et

\begin{proof} Assume that $V=U+ C_{2}(V)$.
Denote by $A$ the vertex Poisson algebra
$\gr_{\E}(V)$ obtained in Theorem \ref{tmain1}.
In particular, $A$ is an $\N$-graded (unital) differential algebra. 
Recall that $A=\coprod_{n\in\N}A_{(n)}$, where $A_{(n)}=E_{n}/E_{n+1}$ for
$n\in \N$.

Let $K$ be the subspace of $V$, spanned 
by those vectors in (\ref{emain-span}). Clearly, $K$ is a graded subspace.
For $m\ge 0$, set $K_{m}=K\cap E_{m}$.

For any linear operator $F$ on a vector space and for any nonnegative
integer $n$, we set $F^{(n)}=F^{n}/n!$. 
{}From Corollary \ref{cdiff}, for any $m\ge 0$, $E_{m}/E_{m+1}$ is 
linearly spanned by the vectors
$$\partial^{(k_{1})}(u^{(1)}+E_{1})\cdots \partial^{(k_{r})}(u^{(r)}+E_{1})$$
for $r\ge 1,\; u^{(i)}\in U,\; k_{1}>k_{2}>\cdots >k_{r}\ge 0$ 
with $k_{1}+\cdots +k_{r}=m$.
By definition we have
\begin{eqnarray*}
\partial^{(k_{1})}(u^{(1)}+E_{1})\cdots \partial^{(k_{r})}(u^{(r)}+E_{1})
&=&(\D^{(k_{1})}u^{(1)}+E_{k_{1}+1})\cdots (\D^{(k_{r})}u^{(r)}+E_{k_{r}+1})\\
&=&u^{(1)}_{-1-k_{1}}\cdots u^{(r)}_{-1-k_{r}}{\bf 1}+E_{m+1}.
\end{eqnarray*}
It follows that $E_{m}=K_{m}+E_{m+1}$. Then
$$V=E_{0}=K_{0}+K_{1}+\cdots +K_{n}+E_{n+1}\subset K+E_{n+1}$$
for any $n\ge 0$.  Since $K$ and $E_{n+1}$ are graded subspaces and
since $E_{m}\subset \coprod_{k\ge 2t+n-1}V_{(k)}$
for $m\ge (n-2)2^{n}$  by (\ref{ekey-need}),
we must have
$$V=K=K_{0}+K_{1}+K_{2}+\cdots,$$
proving the desired spanning property.

Conversely, assume the spanning property.
Notice that if $r\ge 2$, we have $n_{1}\ge 2$, so that
$u^{(1)}_{-n_{1}}\cdots u^{(r)}_{-n_{r}}{\bf 1}\in C_{2}(V)$.
If $n_{r}\ge 2$, we also have
$u^{(1)}_{-n_{1}}\cdots u^{(r)}_{-n_{r}}{\bf 1}\in C_{2}(V)$.
Then we get $V\subset U+C_{2}(V)$, proving $V=U+C_{2}(V)$.
\end{proof}

By slightly modifying the proof of Theorem \ref{tmain-appl}
we immediately obtain the following result (cf. [KL]):

\bt{tpbw}
Let $V=\coprod_{n\ge t}V_{(n)}$ be a lower truncated 
$\Z$-graded vertex algebra such as a vertex operator algebra in the
sense of [FLM] and [FHL] and 
let $S$ be a graded subspace of $V$ such that $\{ u+C_{2}(V)\;|\; u\in S\}$ 
generates $V/C_{2}(V)$ as an algebra. 
Then $V$ is linearly spanned by the vectors
$$u^{(1)}_{-n_{1}}\cdots u^{(r)}_{-n_{r}}{\bf 1}$$
for $r\ge 0,\; u^{(1)},\dots,u^{(r)}\in S,\; n_{1}\ge \cdots \ge n_{r}\ge 1$.
Furthermore, if $S$ is linearly ordered, $V$ is linearly spanned 
by the above vectors with $u^{(i)}>u^{(i+1)}$ when $n_{i}=n_{i+1}$.
\et

\bd{dgenerating-types}
{\em Let $S$ be a subset of a vertex algebra $V$.
We say that $S$ is a {\em type $0$ generating subset of $V$} if 
$V$ is the smallest vertex subalgebra containing $S$, 
$S$ is a {\em type $1$ generating subset of $V$} if  $V$ is linearly
spanned by the vectors
\begin{eqnarray}
u^{(1)}_{-k_{1}}\cdots u^{(r)}_{-k_{r}}{\bf 1}
\end{eqnarray}
for $r\ge 0,\; u^{(i)}\in S,\; k_{i}\ge 1$.
$S$ is called a {\em type $2$ generating subset of $V$} 
if for any linear order on $S$ (if $S$ is a vector space, replace $S$
with a basis), $V$ is linearly spanned by the above vectors with
$u^{(i)}>u^{(i+1)}$ when $n_{i}=n_{i+1}$.}
\ed

\br{rkac}
{\em A type $0$ generating subset is just a generating subset
in the usual sense and a type $1$ generating subset of $V$ is also called
a strong generating subset $V$ in [K].}
\er

\bt{tstrong}
Let $V$ be a lower truncated $\Z$-graded vertex algebra 
and let $U$ be a graded subspace. Then the following three statements
are equivalent: (a)
$U$ is a type $1$ generating
subspace of $V$. (b) $U$ is a type $2$ generating
subspace of $V$. (c)
$U/C_{2}(V)=\{u+C_{2}(V)\;|\; u\in U\}$ 
generates $V/C_{2}(V)$ as an algebra.
\et

\begin{proof} By definition, (b) implies (a) and 
by Theorem \ref{tpbw}, (c) implies both (a) and (b).
Now it suffices to prove that (a) implies (c).
Assuming (a) we have that $V/C_{2}(V)$ is linearly spanned by
the vectors $u^{(1)}_{-k_{1}}\cdots u^{(r)}_{-k_{r}}{\bf 1}+C_{2}(V)$
for $r\ge 0,\; u^{(i)}\in U,\; k_{i}\ge 1$.
If $k_{i}\ge 2$ for some $i$, we have
$u^{(1)}_{-k_{1}}\cdots u^{(r)}_{-k_{r}}{\bf 1}\in C_{2}(V)$. Then
$V/C_{2}(V)$ is linearly spanned by
the vectors $u^{(1)}_{-1}\cdots u^{(r)}_{-1}{\bf 1}+C_{2}(V)$
for $r\ge 0,\; u^{(i)}\in U$. That is,
$U/C_{2}(V)$ generates $V/C_{2}(V)$ 
as an algebra.
\end{proof}


With Lemma \ref{lcheck-assumption}, from the proof of Theorem \ref{tmain-appl} 
we immediately have:

\bp{pmodule} 
Let $V$ be a lower truncated $\Z$-graded vertex algebra
and let $U$ be a graded subspace of $V$ such that $U$ generates $V/C_{2}(V)$
as an algebra. Let $W$ be a lower truncated $\Z$-graded $V$-module and let $W^{0}$ be a graded subspace
of $W$ such that $W=W^{0}+C_{2}(W)$.  Then $W$ is spanned by the
vectors
\begin{eqnarray*}
u^{(1)}_{-1-k_{1}}\cdots u^{(r)}_{-1-k_{r}}w
\end{eqnarray*}
for $r\ge 1,\; u^{(1)},\dots, u^{(r)}\in U,\; 
w\in W^{0},\;k_{1}> \cdots > k_{r}\ge 0$.
\ep

\end{document}